\documentclass[12pt, Letterpaper]{article}
\usepackage[utf8]{inputenc}
\usepackage[english]{babel}
\usepackage[nottoc]{tocbibind}
\usepackage{graphicx,caption}
\usepackage{amssymb, amsmath, mathtools, amsfonts}
\usepackage{longtable, supertabular}
\usepackage{url, rotating, array, multirow, geometry}
\usepackage{float}
\usepackage{ytableau, youngtab}
\usepackage{subfloat}
\usepackage{sidecap}
\usepackage{ragged2e}
\usepackage{mdwlist}
\usepackage{enumitem}
\usepackage{setspace}
\usepackage{hyperref}
\usepackage{color}
\usepackage{subfig}
\usepackage[labelfont=bf]{caption}
\usepackage[compact]{titlesec}
\begin{document}\hypersetup{pdfborder={000}}
\pagenumbering{gobble}
\title{A Note on Frobenius Character Formula}
\author{
Inayat R. Bhat\footnote{G. D. C. Shopian, inayatrb@gmail.com}
}
\date{\today}
\maketitle
\thispagestyle{plain}
\begin{abstract}
\noindent In this note we give an alternative derivation of Frobenius Character Formula using wedge product, the resulting formula is used to calculate the character table of $S_{2}$ and $S_{3}$ Group which reproduces their character table correctly, calculations for classes of $S_{3}$ Group are given in the appendix \ref{app:02}. The calculations are given in terms Schur Functions. In the first section, a brief indicative review of Symmetric Group and related formulas is given. In Section second we derive the formula. A naive approach to computational complexity is also presented at the end.
\end{abstract}
\newpage
\listoftables
\tableofcontents
\newpage
\section{Introduction}
\noindent A Symmetric group $Sym(X)$ or $S_{n}$ on a finite set $X$ of $n$ objects is a collection of all bijections (permutations) of set $X$ onto itself . We introduce some notation and point to some general formulas about $S_{n}$ in a less formal way. Detailed subsequent material can be found in \cite{ref:B7}, \cite{ref:B12}, \cite{ref:B4},\cite{ref:B2},\cite{ref:B3}     - \cite{ref:B6} \\
For a finite set of arbitrary elements of cardinality $n$,  $X = \{ x_{k} : k \leq n \in \mathbb{N}\} \equiv \{ k : k \leq n \in \mathbb{N}\} $,  the set of all bijections $S_{n} : X \rightarrow X$ has cardinality $n!$, thus, $ \forall k$, the only identity bijection, which maps each element to itself is $I_{d}: k \rightarrow k$ and can be denoted by $(1)(2)(3)...(n)$, which is an identity permutation, where each parenthesis $(i)$ symbolises that  $i$ is mapped to $i$. For brevity this is denoted by $(1^{n})$ where each number is replaced by $(1)$ depicting a cycle of length one with each element mapped to itself. A bijection which maps a pair $(j,k) \rightarrow (k,j)$ is called a transposition. One single transposition $(jk)$, which is a permutation of symbols $j, k$ is denoted by $(2)$ (meaning mapped to each other) and there are $\dfrac{n(n-1)}{2}$ such bijections and a pair of transpositions is denoted by $(ij)(lm) \equiv (2^{2})$, where $i, j, l, m \in X$. Thus a bijection having $\nu_{2}$ transpositions and the rest mapped to themselves is denoted by $(1^{\nu_1} 2^{\nu_{2}})$ with $\nu_{1} + 2\nu_{2} =n$. With this process a generic structure of a bijection in this cyclic notation is $(l) = (1^{\nu_{1}} 2^{\nu_{2}} 3^{\nu_{3}}\cdots k^{\nu_{k}})$ with $\sum_{k=1}^k k \nu_{k}   = n$, where $\nu_{k} \in \mathbb{N}$.  If we define an equivalence relation $ \forall k, \sigma_{k} \in S_{n}, \sigma_{k} \sigma_{i} \sigma_{k}^{-1} = \sigma_{j}$, where for instance $\sigma_{k} = (ijk...lmn)$, then $\sigma_{k}^{-1} = (nml...kji)$, the total number of elements in $S_{n}$ fall into disjoint sets or equivalence classes with each element in a class having the same cyclic structure and number of elements in each class equals the number as computed from formula \eqref{eq:1}. These classes are often called conjugacy classes.\\
For a generic class, define a set 
\begin{equation}
K= \{1,2,...,k\} \label {s:1}
\end{equation}
 where elements in $K$ are cycle types available in class $K \subset \{k: 1\leq k \leq |X| =n\}$. It is desirable for the purpose here to write the elements of $K$ in increasing or decreasing sequence and $K^{c}$ is the complement of set $K$ with $K\cup K^{c} =X$. The number of permutations of a particular cycle structure $(1^{\nu_1} 2^{\nu_{2}} 3^{\nu_{3}}...k^{\nu_{k}})$ is given by 
\begin{equation}
N_{(l)} = \dfrac{n!}{\prod_{r\in K} r^{\nu_{r}}\nu_{r}!} \label {eq:1}
\end{equation} \\
Each cycle of length $r$ can be written in $r$ different ways, therefore division by $r$ cancels repetitions, for $\nu_{r}$ such cycles we divide by $r^{\nu_{r}}$, further, the $\nu_{r}$ cycle can be permuted in $\nu_{r}!$ ways. Formula \ref{eq:1} is used to compute the number of elements in a given conjugacy class. Scaling $k$ to unity in $\sum_{k =1}^n k \nu_{k} = n$ as $\sum_{k=1}^n  \nu_{k} + \sum_{k=2}^n  \nu_{k} + \sum_{k=3}^n  \nu_{k} + ... + \sum_{k=j}^n  \nu_{k} + ... +  \sum_{k=n}^n \nu_{n} = n$, and write 
\begin{equation}
\sum_{j=1}^n \lambda_{j}= n \label{eq:2}
\end{equation}
 with 
\begin{equation}
\lambda_{j} = \sum_{k=j}^n  \nu_{k}  \label{eq:3}
\end{equation}
The set of numbers $(\lambda) = (\lambda_{1}, \lambda_{2}, \lambda_{3}, ..., \lambda_{j}, ..., \lambda_{n})$ is partition of $n$ denoted as $\lambda \vdash n$, From equation \ref{eq:3} since, $\sum_{k=j}^n  \nu_{k} \geq \sum_{k=j+1}^n  \nu_{k}$, we have 
\begin{equation}
\forall j, \quad \lambda_{j}\geq \lambda_{j+1} \label{eq:4}
\end{equation}
Given a partition $(\lambda)$ then $\nu_{k} = \lambda_{k} - \lambda_{k+1}$. Therefore, From a partition $(\lambda)$ of $n$ and the cycle structure of $S_{n}$ can be recovered, and hence have a one to one correspondence.The total number of classes of $S_{n}$ is in one to one correspondence with the integer partitions of $n$, which are solutions of equation (\ref{eq:2}) subjected to condition (\ref{eq:4}). Upto now we represented the bijections $S_{n} : X \rightarrow X$ as permutations in cyclic notation.\cite{ref:B2},\cite{ref:B3}\\
It is desirable to represent the elements of  $S_{n}$ by a set of $d \times d $ invertible matrices $M_{ij}$, and can be generally taken as unitary matrices, where $1 \leq d \leq n!$, for example an identity can be represented by a matrix $\delta_{ij}$ where $\delta_{ij}$ is Kronecker delta, these matrices can be often reduced for arbitrary $d$ to matrices of block diagonal form using similarity transformation, these blocks of matrices with dimensions $\leq d$ cannot be reduced further into sub-blocks, these blocks are called irreducible representations. Using Schur Lemma \cite{ref:B2}, \cite{ref:B4},\cite{ref:B3}, if $\hat{M}_{ij}$ is a transformed matrix, then $\forall \sigma_{k} \in S_{n}, \hat{M}_{ij} (D_{ii} - D_{jj}) =0$, where $D_{ii}$, $D_{jj}$ are diagonal and the values at diagonal points of these matrices determine the Block diagonal form of $\hat{M}_{ij}$. if $\Gamma$ is reducible then in terms of irreducible representations, $\Gamma = \bigoplus_{k} m_{k} \Gamma_{k}$, where $m_{k}$ is the multiplicity of a $\Gamma_{k}$, and gives the number of times a block matrix is repeated in the process. Further this direct sum is unique apart from the arrangement of blocks.\\
There is a one to one correspondence between the number of irreducible representations of a group and the total number of classes. The number of irreducible representations for a group equals the number of its classes.\\
The representation of group can be realised in two equivalent ways, first is group homomorphism which is a map that preserves group multiplication. if h is a homomorphism from $(G, .)$ to $(H, \times)$, then $\forall g_{i}, g_{j} \in$ $G$, $h(g_{i} . g_{j}) = h(g_{i}) \times h(g_{j})$. A representation of $S_{n}$ on a vector space $V$ of dimensions $d$ is a homomorphism $\rho_{\sigma}:S_{n}\rightarrow GL(d, V)$, where $\sigma \in S_{n}$ and $d$ are dimensions of $V$, which is choosen as the underlying vector space, as for example the $V$ can be the group itself say $S_{n}$ then choosing its elements as basis of the underlying vector space we get $d = n!$, representations with $d > n!$ are possible but with the choice $d = n!$ is the minimum dimension which contains all the irreducible representations and is called regular representation.  Second, any representation defines an action of the group on vector space $V$ which is a linear map $\Phi : S_{n} \times V \rightarrow V$, called left action if $\forall g, g_{1},g_{2} \in S_{n}, \|x\rangle, \|y\rangle \in V, \alpha, \beta \in \mathbb{C}$, 
\begin{enumerate}
\item $\Phi (e, |x\rangle) = |x\rangle$
\item $ \Phi( g_{1},g_{2}, |x\rangle) = \Phi (g_{1}, \Phi(g_{2},|x\rangle))$
\item $\Phi(g, \alpha |x\rangle + \beta |y\rangle) = \alpha \Phi(g, |x\rangle) + \beta \Phi(g, |y\rangle)$
\end{enumerate}
the $\Phi$ map, therefore is compatible with group multiplications and operations in $V$. The right action $\Psi$, is defined similarly and in particular $\Psi(|x\rangle, g)= \Phi(g^{-1}, |x\rangle)$. For instance, the left action $\Phi: S_{n} \times $B$ \rightarrow GL(\mathbb{R}^{m})$ is, where $B = \{| e_{k}\rangle\}_{k=1}^{m}$ is the standard basis set of $\mathbb{R}^{m}$.
\begin{equation}
\Phi_{\sigma_{i}} \left |e_{k}\right\rangle = \sum_{j=1}^{m} \left |e_{j}\right\rangle \left\langle e_{j} \right | \Phi_{\sigma_{i}} \left |e_{k}\right\rangle = \sum_{j=1}^{m} \left |e_{j}\right\rangle \left\langle e_{j} \right | e_{k} . \sigma_{i} \rangle
\end{equation}
 where$\left\langle e_{j} \right | e_{k} . \sigma_{i} \rangle = \left \langle e_{j} \right | e_{\sigma_{i}(e_{k})} \rangle = \delta_{j,\sigma_{i}(e_{k})}$ are the matrix elements, where Dirac ket-bra notation is used.\\
 if $M_{ij}$ denotes a matrix of a representation of an element in $S_{n}$ then the Map $\chi: S_{n} \rightarrow \mathbb{C}$ is called character and is equal to $\chi=\sum_{i} M_{ii} $, if $M_{ij}$ represents a matrix of irreducible representation the $\chi$ is called simple, a compound character is a sum of simple characters as reducible representation is a sum of irreducible representations.\\
A table of $\chi^{(\lambda)}_{(l)}$ with columns (rows) given by partitions $(\lambda)$ and rows (columns) by class structure $(l)$ is called character table. \cite{ref:B3}\\
The partitions of an integer $\lambda \vdash n$ can be viewed geometrically as a 2-D array of boxes which are left justified and each row has $\lambda_{i}$ boxes called Ferrers diagram or Young Tableau or diagram. A Ferrers diagram is filled with dots and a Young diagram by numbers. if $\lambda$ is a tableau of certain shape, then the coordinates of a box $(i,j)$ is its row and column number with top left set to $(1,1)$, along the rows $i$ increases by 1, down the column $j$ increases by 1. A tableau with $n$ boxes can be filled with a set of numbers according to the two different rules which in turn give two types of Young tableau as:
 \begin{enumerate}
\item Standard Young Tableau (SYT): if the set of numbers chosen to fill a tableau of $n$ boxes are distinct and filled $s.t.$ $\forall i \in \lambda_{i}, (i,j) < (i,j+1)$ and $\forall j \in \lambda_{j}, (i,j) < (i +1,j)$.
\item Semi-Standard Young Tableau (SSYT): if the boxes are filled with numbers $s.t.$ $\forall i \in \lambda_{i}, (i,j) \leq (i,j+1)$ and $\forall j \in \lambda_{j}, (i,j) < (i+1, j)$.
 \end{enumerate}
 Therefore, $1.$ is strictly increasing across rows and down columns and $2.$ is weakly increasing across rows and strictly increasing down columns.\\
For a given shape $\lambda \vdash n$, the total number of SYT that can be obtained equals the dimensions of the irreducible representation corresponding to $\lambda$, and is given by hook length formula $f^{(\lambda)} = \dfrac{n!} { \prod_{{(i,j}\in \lambda} h(i,j)}$, where $h(i,j)$ is the hook associated with box $(i,j)$. The hook of a box $(i,j)$ is the number of boxes obtained by adding together the box to number of boxes to its right and number of boxes below it.\\
For the SSYT, since they allow repetitions of symbols along rows, they turn out to be useful, with each symbol $k$ in the tableau associate a variable $x_{k}$, and a monomial $x^{T} = x_{1}^{w_{1}}x_{2}^{w_{2}}...x_{k}^{w_{k}}$, where the set $w(T) = \{ {w_{1}},{w_{2}}, ..., {w_{k}} \},$ is called weight of a SSYT. The weight $w_{k}$ is simply the number of times the integer $k$ appears in the SSYT.\\
 A symmetric function is a function which is invariant under the permutation of variables. If $\lambda$ is a partition or a shape of a tableau, then the symmetric function of type \cite{ref:B7}
 \begin{equation}  
 s_{(\lambda)} = \sum_{T} x^{w(T)} \label{eq:5}
\end{equation}
 is called a Schur function. the set of all Schur functions forms a basis for a ring of symmetric functions $\Lambda^{n}$ if $\lambda \vdash n$. Other type of basis is power basis labeled by $(l)$, the cycle structure of a class, where,
\begin{equation}
p_{r} = \sum_{k=1}^{n} x_{k}^{r}, \quad r\in K,  \label{eq:6}
\end{equation} 
These are a set of linearly independent quantities, the multinomial expansion of $p_{r}^{\nu_{r}}$ is given by \cite{ref:B4}
\begin{equation}
p_{r}^{\nu_{r}} = \displaystyle\sum_{\sum_{j} k_{r,j} = \nu_{r}} \dfrac{\nu_{r}!}{\prod_{j} k_{r,j}!} \prod_{j} x_{j}^{rk_{r,j}} \label{eq:7}
\end{equation}
and find the product of these symmetric polynomials for a generic class $(l)$, $p_{(l)} = \prod_{r} p_{r}^{\nu_{r}}$, which is given by \cite{ref:B4}
\begin{equation}
p_{(l)} = \prod_{r} p_{r}^{\nu_{r}} =\sum_{\tiny\begin{gathered}\forall r, \sum_{j} k_{r,j} = \nu_{r},\\ \forall j, \sum_{r} rk_{r,j} =\lambda_{j}\end{gathered}} \prod_{r\in K} \dfrac{\nu_{r}!}{\prod_{j} k_{r,j}!} \prod_{j} x_{j}^{\lambda_{j}} \label{eq:8}
\end{equation}
The coefficients on the right hand side are compound characters \cite{ref:B12}, \cite{ref:B4} with the antisymmetric polynomial defined by the determinant 
\begin{equation}
D(x_{i}) = \prod_{i<j} (x_{i} - x_{j}) \label{eq:9}
\end{equation}
Then the simple characters are coefficients of terms $\prod_{j} x^{\lambda_{j} + n - j}_{\sigma_{(j)}}$ in $p_{(l)}D(x_{i})$, that is
\begin{equation}
p_{(l)} D(x_{i}) = \sum_{(\lambda)} \chi^{(\lambda)}_{(l)} \sum_{\sigma \in S_{n}} \operatorname{sgn}(\sigma) \prod _{j =1}^{n} x^{\lambda_{j} +n -j}_{\sigma(j)} \label{eq:10}
\end{equation}
Formula \ref{eq:10} is Frobenius Formula for characters, the proof is given in \cite{ref:B12}, \cite{ref:B4}, \cite{ref:B6}, The Formula \ref{eq:10} can be alternatively written with the Vandermonde Determinant defined as 
\begin{equation}
V(x_{i}) = \displaystyle \prod_{i<j} (x_{j} - x_{i}) \label{eq:11}
\end{equation}
The Vandermonde equals its transpose $V(x_{i}) = \widetilde{V}(x_{i})$and is related to the determinant $D(x_{i})$ as
\begin{equation}
D(x_{i}) = (-1) ^\frac{n(n-1)}{2} V(x_{i}) \label{eq:12}
\end{equation}
Since each flip of factor $(x_{j} - x_{i}) \rightarrow (x_{i} - x_{j})$ costs a minus sign there are $\binom{n}{2}$ flips.\\
Further, for a sequence $\{\alpha_{k}\}_{k =1}^{n}$, the determinant of the form denoted by $a_{\alpha_{1}, \alpha_{2} ,\cdots, \alpha_{n}}$ is called an alternant \cite{ref:p1}, \cite{ref:B6}  and equals
\[
\begin{vmatrix}
  x_{1}^{\alpha_{1}} &  x_{1}^{\alpha_{2}}  & \cdots  &x_{1}^{\alpha_{n}} \\
  x_{2}^{\alpha_{1}} &  x_{2}^{\alpha_{2}}  & \cdots  &x_{2}^{\alpha_{n}} \\
     \vdots & \vdots &\vdots\\
    x_{n}^{\alpha_{1}} &  x_{n}^{\alpha_{2}}  & \cdots  &x_{n}^{\alpha_{n}} \\
\end{vmatrix}
\]
For a special case where the sequence $\{\alpha_{k}\}_{k=1}^{n}$ is $\{ (n-k)\}_{k=1}^{n}$, which is a stair case often denoted by $\delta = (n-1,n-2,\cdots,1,0)$ then the alternant $a_{\delta}$ , Vendermonde $V(x_{i})$ are relates as
\begin{equation}
a_{\delta} = (-1)^{\lfloor \frac{n}{2} \rfloor} \widetilde{V}(x_{i})=(-1)^{\lfloor \frac{n}{2} \rfloor} {V}(x_{i})
\end{equation}
where $\lfloor \dfrac{n}{2} \rfloor$ is the floor of $\dfrac{n}{2}$. For a partition $\lambda \vdash n$, the alternant $a_{\lambda + \delta}$ is defined similarly and equals  equation \ref{eq:10}. Further Schur functions as in \ref{eq:5} are related to alternants $a_{\lambda +\delta}$ and $a_{\delta}$ by bi-alternant formula \cite{ref:p1},\cite{ref:B7}, 
\begin{equation}
s_{(\lambda)} = \dfrac{a_{\lambda+\delta}}{a_{\delta}} \label{eq:14}
\end{equation}
The transformation from power sum basis as given in equation \ref{eq:8} and Schur basis as in equation \ref{eq:5} and \ref{eq:14} is given by
\begin{equation}
p_{(l)} = \sum_{\lambda \vdash n} \chi^{(\lambda)}_{(l)} s_{(\lambda)} \label{eq:15}
\end{equation}
where the sum is over all partitions and $ \chi^{(\lambda)}_{(l)}$ are entries in the character table along the colomn(row) for constant $(l)\in K$,  see equation \ref{s:1} and its inverse transformation is given by \cite{ref:B7}
\begin{equation}
s_{(\lambda)} = \sum_{l \vdash n} \dfrac{\chi^{(\lambda)}_{(l)}}{\prod_{r\in K} r^{\nu_{r}}\nu_{r}! } p_{(l)} \label{eq:16}
\end{equation}
where $\prod_{r\in K} r^{\nu_{r}}\nu_{r}! $ is the normalisation factor and from equation \ref{eq:15} and \ref{eq:16}, the character table $\chi^{\lambda}_{(l)}$ is the transition matrix between power basis $p_{(l)}$ and  Schur basis $s_{(\lambda)}$. 
\section{Methedology}
\noindent Consider equation \ref{eq:7} and find $d(p_{r}^{\nu_{r}}) = \dfrac{\nu_{r}}{p_{r}} p_{r}^{\nu_{r}}$, where $\nu_{r}$ is the exponent of cycle $r$ in the generic cycle structure of class $(1^{\nu_1} 2^{\nu_{2}} 3^{\nu_{3}}...k^{\nu_{k}})$, thus $r \in K$ and this equals
\begin{equation}
\nu_{r} p_{r}^{\nu_{r}-1} dp_{r} = \displaystyle\sum_{\sum_{j} k_{r,j} = \nu_{r}} \dfrac{\nu_{r}!}{\prod_{j} k_{r,j}!} \prod_{j} x_{j}^{rk_{r,j}} \sum_{j=1}^{n} \dfrac{r k_{r,j}}{x_{j}} 	\label{eq:17}
\end{equation} 
multiplying by $\dfrac{p_{r}}{\nu_{r}}$ and taking the $\bigwedge$ product \cite{ref:B11} of all the cycle types in the class for $r \in K$,
we get on one side
\begin{equation}
\bigwedge_{r \in K} p_{r}^{\nu_{r}} dp_{r} = \prod_{r\in K} p_{r}^{\nu_{r}} \bigwedge_{r\in K} dp_{r} \label{eq:18}
\end{equation}
On the other side of equation \ref{eq:7}, the $\bigwedge$ product is coefficient times $det \left(\dfrac{rk_{r,j}}{x_{j}}\right) \bigwedge_{r\in K} dx_{r}$ where $det\left(\dfrac{rk_{r,j}}{x_{j}}\right)$ is a determinant of dimensions $k\times n$, $k$ is cardinality of $K$,
\begin{equation}
\prod_{r\in K}p_{r}^{\nu_{r}} \bigwedge_{r\in K} dp_{r} = \prod_{r\in K} p_{r}\displaystyle\sum_{ \forall r \in K, \sum_{j} k_{r,j} = \nu_{r}} \prod_{r \in K}\left(\dfrac{(\nu_{r} -1)!}{\prod_{j=1}^n k_{r,j}!}\right) \prod_{j=1}^{n} x_{j}^{\sum_{r=1}^{n} rk_{r,j}-1}  d(k_{r,j}) \bigwedge_{r\in K} dx_{r}  \label{eq:28-1}
\end{equation}
for those $r \in K^{c}$, use equation \ref{eq:6}, and find $dp_{r} = r \sum_{j =1}^{n}x_{j}^{r-1}dx_{j}$, and take wedge with those $r \in K$ in \ref{eq:28-1}
\begin{equation}
\prod_{r\in K} p_{r}^{\nu_{r}} \bigwedge_{r=1}^{n} dp_{r} =\prod_{r\in K}p_{r}\displaystyle\sum_{ \forall r \in K, \sum_{j} k_{r,j} = \nu_{r}} \prod_{r \in K}\left(\dfrac{(\nu_{r} -1)!}{\prod_{j=1}^{n} k_{r,j}!}\right) \prod_{j=1}^{n} x_{j}^{\sum_{r=1}^{n} rk_{r,j}}  d(k_{r,j},x_{j}) \bigwedge_{r =1}^{ n} dx_{r}  \label{eq:19}
\end{equation}
where we define a determinant $d(k_{r,j},x_{j})$ as 
\[
d(k_{r,j},x_{j}) = 
\begin{cases}
\dfrac{rk_{r,j}}{x_{j}} &  r \in K\\
r x_{j}^{r-1} &  r \in K^{c}
\end{cases} 
\]
For the product $\bigwedge_{r=1}^{n} dp_{r}$ on LHS of \ref{eq:19}, we again use \ref{eq:6} and find 
\begin{equation}
\bigwedge_{r=1}^{n} dp_{r} = \mathcal{J} \bigwedge_{r=1}^{n} dx_{r} \label{ref:eq20}
\end{equation}
where $\mathcal{J}$ is the Jacobian determinant $\mathcal{J} = \dfrac{\partial p_{r}}{\partial x_{j}} = n!V(x_{j})$. With these inclusions in \ref{eq:19} and $ p_{(l)} =\prod_{r \in K} p_{r}^{\nu_{r}}$ and $p_{(\lambda)} = \prod_{r\in K}p_{r}$, and cancelling the n! form both sides (n! from the right side comes from det $d(k_{r,j},x_{j}) ) $and $\bigwedge_{r=1}^{n} dx_{r}$, we get 
 \begin{equation}
p_{(l)} V(x_{j}) = p_{(\lambda)} \displaystyle\sum_{ \forall r \in K, \sum_{j} k_{r,j} = \nu_{r}} \prod_{r \in K}\left(\dfrac{(\nu_{r} -1)!}{\prod_{j=1}^{n} k_{r,j}!}\right) \prod_{j=1}^{n} x_{j}^{\sum_{r=1}^{n} rk_{r,j}} sgn(\sigma)d(k_{r,j},x_{j}) \label{eq:21}
\end{equation}
where $sgn(\sigma)$ is the sign of the $det \left(d(k_{r,j},x_{j})\right)$ when it is arranged into a two block form, an upper block of constants  as $\dfrac{ k_{r,j}}{x_{j}},  r \in K$ and a lower block of variables as $x_{j}^{r-1} ,  r \in K^{c}$.\\
$p_{(\lambda)} = \prod_{r\in K}p_{r} = \prod_{r=1}^{l(\lambda)} \sum_{k} x_{k}^{\lambda_{k}}$, where $l(\lambda)$ is partition length for $\lambda = (\lambda_{1},\lambda_{2}, \cdots \lambda_{|K|})$, as $l(\lambda) = |K|$, when elements in $K$ are arranged in decreasing sequence, and $\sum_{r=1}^{n} r k_{r,k}$ is not set equal to $\lambda_{r}$ because the equation employs these internal components $k_{r,j}$ to compute the polynomial arising from the sum in equation \ref{eq:21} ${ \forall r \in K, \sum_{j} k_{r,j}} = \nu_{r}$.\\
Consider now the sign $sgn(\sigma)$ of the $det (k_{r,j},x_{j})$ as follows. Consider the underlying set $X$ of $S_{n}$ and $K\subset X$,  $k \in X$  with  $1 \leq k \leq n \in \mathbb{N}$, Further $\sum_{k =1}^{|K|} \nu_{k} k = n$, $\nu_{k}$ is the frequency of a cycle of type $k$ placing the elements of the set in an increasing sequence, by swapping two at a time and if each swap costs a negative sign, then the sign in $sgn \sigma = (-1)^{n(t)}$, where $n(t)$ is the number of transpositions in the cyclic permutation of elements in $K$ indexed with elements of $X$. As an example if $K = \{2,4,6\}$, means if the cycle structure is of the form $(l) = (2^{\nu_{2}} 4^{\nu_{4}} 6^{\nu_{6}})$ then indexing $K$ with elements $X$, $(124)(365)$ is the required permutation with $4$ transpositions.\\
Equation \ref{eq:21} can be compared with equation \ref{eq:10} apart from the sign which arises from using equation \ref{eq:12}, and factoring the $V(x_{j})$, verifies equation \ref{eq:15}\\
The illustration of the method is given below for simple cases of $S_{2}$. Details for $S_{3}$ are given in Appendix \ref{app:02} for all three classes. \\
For $S_{2}$, first choose the identity class $(l) = (1^{2})$, $p_{(1^{2})} = p_{1}^{2} = (x_{1} + x_{2})^{2}$, therefore, expanding the $p_{1}^{2} =(x_{1} + x_{2})^{2} = \displaystyle\sum_{\sum_{j} k_{1,j} =2} \dfrac{2!}{k_{1,1}! k_{1,2}!} x_{1}^{k_{1,1}}x_{2}^{k_{1,2}}$ and therefore 
\begin{equation}
2p_{1} dp_{1} = \displaystyle\sum_{\sum_{j} k_{1,j} =2} \dfrac{2!}{k_{1,1}! k_{1,2}!} x_{1}^{k_{1,1}}x_{2} \left( \dfrac{k_{1,1}}{x_{1}} dx_{1} + \dfrac{k_{1,2}}{x_{2}}dx_{2}\right) \label{eq:ex1}
\end{equation}
for $p_{2} = (x_{1}^{2} + x_{2}^{2})$, we have $dp_{2} = 2(x_{1}dx_{1} + x_{2}dx_{2})$
\begin{equation}
2p_{1}^{2}dp_{1} \wedge dp_{2}= p_{1} \displaystyle\sum_{\sum_{j} k_{1,j} =2} \dfrac{2!}{k_{1,1}! k_{1,2}!} x_{1}^{k_{1,1}}x_{2}^{k_{1,2}} 
\begin{vmatrix}
  \dfrac{k_{1,1}}{x_{1}} &  \dfrac{k_{1,2}}{x_{2}}\\
  2x_{1} &  2x_{2}  \\
\end{vmatrix}
dx_{1} \wedge dx_{2}
\end{equation}
using $dp_{1} \wedge dp_{2} = \begin{vmatrix}
  1 &  1\\
  2x_{1} &  2x_{2}  \\
\end{vmatrix} dx_{1} \wedge dx_{2}$, we get,
\begin{equation}
p_{(1^{2})} (x_{2} -x_{1}) = p_{(1)} \displaystyle\sum_{\sum_{j} k_{1,j} =2} \dfrac{2!}{k_{1,1}! k_{1,2}!} x_{1}^{k_{1,1}}x_{2}^{k_{1,2}} 
\begin{vmatrix}
  \dfrac{k_{1,1}}{x_{1}} &  \dfrac{k_{1,2}}{x_{2}}\\
  2x_{1} &  2x_{2}  \\
\end{vmatrix} \label{eq:24}
\end{equation}
Therefore, from equation \ref{eq:25} and Table \ref{tab:t1} we have
\begin{equation}
p_{(1^{2})} (x_{2}-x_{1}) = p_{(1)} (x_{2}-x_{1})(x_{2}+x_{1}) \label{eq:25}
\end{equation} 
where $p_{(1)} = (x_{1} +x_{2})$, if we cancel the common factor we get $p_{(1^{2})} = x_{1}^2 +  x_{2}^2 +  2x_{1} x_{2}$ and thus 
\begin{equation}
p_{(1^{2})} = s_{(2)} + s_{(1,1)} \label{eq:26}
\end{equation}
\begin{longtable}{p{0.10\textwidth}|p{0.30\textwidth}}
\caption{Table for $(1^{2})$ in $S_{2}$}
\label{tab:t1}\\
\hline
  ($k_{1,1}, k_{1,2})$ &$p_{(1^{2})} (x_{2} -x_{1}),$  \ref{eq:24} \\
\hline
 (2,0)&$p_{(1)} x_{1}x_{2}$			   \\
 (0,2)&$-p_{(1)} x_{1}x_{2}$			   \\
 (1,1)&$p_{(1)} (x_{2}^{2} - x_{1}^{2})$	    \\
\hline
&$p_{(1)} (x_{2}^{2} - x_{1}^{2})$	\\	
\hline			
\end{longtable}
\noindent First few schur functions are listed in Table \ref{t:a1} of Appendix \ref{app:01}. Schur functions are stable \cite{ref:p1}, we simply put $x_{3} =0$ to arrive at $s_{(\lambda)}(x_{1},x_{2})$ from $s_{(\lambda)}(x_{1},x_{2},x_{3})$. If, however in \ref{eq:25}, we use $D(x_{1},x_{2}) = (x_{1} -x_{2})$ as in \ref{eq:9}, \ref{eq:10}, the polynomial expansion is \begin{center}$x_{1}^{3} + x_{1}^{2}x_{2} -x_{1}x_{2}^{2} - x_{2}^{3}$\end{center}Then for $n=2$, $\delta = (1,0)$, $(\lambda) = (2,0)$, the coefficient of $a_{\lambda + \delta} = a_{(3,0)}$ is $1$ in $x^{3}$ term. Similarly for the other class of $S_{2}$, which is $(2)$, we arrive at $p_{(2)} = x_{1}^{2} + x_{2}^{2}  = x_{1}^{2} + x_{2}^{2} + x_{1}x_{2} - x_{1}x_{2} $
\begin{equation}
p_{(2)} = s_{(2)} -s_{1,1} \label{eq:27}
\end{equation}    \\
The polynomial expansion is \begin{center} $x_{1}^{3} - x_{1}^{2}x_{2} +x_{1}x_{2}^{2} - x_{2}^{3}$\end{center} Equations \ref{eq:26} and \ref{eq:27} verify \ref{eq:15}. Calculations for $S_{3}$ classes are given in \ref{app:02} see also \ref{c:ps}, which are all alternating polynomials and hence are divisible by $V(x_{i})$ and therefore the quotient will be a symmetric polynomial which can be written as a linear sum in Schur basis. The corresponding factors of power basis are kept as such in tables given in \ref{app:02}.\\
The formula \ref{eq:21} can be applied to subgroups, this can be readily seen from equation \ref{eq:17}, which has a reduced power for $r^{th}$ cycle by $1$, if we refrain from multiplying by $p_{r}$, we get a class which has one less such cycle. Similarly, the total reduction will be by $\sum_{r}r \in K$. The generic class structure will be $(l') = (1^{\nu_{1}-1} 2^{\nu_{2}-1} 3^{\nu_{3}-1}\cdots k^{\nu_{k}-1})$, with $(l') \in S_{n - \frac{k(k+1)}{2}}$. Further, if we donot include $\wedge$ products from $\bigwedge_{r \in K^{c}}dp_{r}$,
while arriving at \ref{eq:19}, we get an equation 
\begin{equation}
\prod_{r\in K}p_{r}^{\nu_{r}-1} \bigwedge_{r\in K} dp_{r} = \displaystyle\sum_{ \forall r \in K, \sum_{j} k_{r,j} = \nu_{r}} \prod_{r \in K}\left(\dfrac{(\nu_{r} -1)!}{\prod_{j=1}^n k_{r,j}!}\right) \prod_{j=1}^{n} x_{j}^{\sum_{r=1}^{n} rk_{r,j}-1}  d(k_{r,j}) \bigwedge_{r\in K} dx_{r}  \label{eq:28}
\end{equation}
To this equation \ref{eq:28}, multiply by $\bigwedge_{r \in K^{c}} dp_{r}$ from equation \ref{eq:6} but for $r \in K^{c}$. we arrive at the formula for character of subgroup of $S_{n}$ of reduced power. Therefore, for subgroups we use \ref{eq:19} without multiplication by $p_{(\lambda)} = p_{(k,\cdots3,2,1)}$ Applying this to class $(1^{2})$ of $S_{2}$ as in Table \ref{tab:t1} and without multiplying by factor $p_{1}$, we get, following procedure from equation \ref{eq:ex1}  $p_{1} \in S_{1} = x_{1} + x_{2} = s_{(1)}$, which is simply a Young diagram of one box corresponding to partition of $(1)$, The other class $(2)$ of $S_{2}$ gives $s_{\phi}=1$ because when we differentiate $p_{2}$ with one power we are left with $1$, The class $(2) \in S_{2}$ and is not contained in $S_{1}$. Similarly $S_{2}$ as a subgroup of $S_{3}$ can be seen from the calculation given at the bottom of each Table \ref{tab:t2}, \ref{tab:t3}, \ref{tab:t4} in terms of $s_{(\lambda)}$ and determinant $D(x_{i})$ of Appendix \ref{app:02}
\section{Computational Complexity}
\noindent In a naive approach, we proceed as follows: The generic class $(l) = (1^{\nu_{1}} 2^{\nu_{2}} 3^{\nu_{3}}\cdots k^{\nu_{k}})$ has a constraint $\sum_{k=1}^k k \nu_{k}   = n$, which is a diaphontine equation but with solutions restricted to $\nu_{k} \in \mathbb{N}_{0}$. Assume, its solutions are given. Therefore, the solution set is all the classes of $S_{n}$ as $\nu_{k}$$'s$ correspond to the cycle type.  These solutions are exactly the sum conditions over rows of determinant in equation \ref{eq:19}, which is  $\forall r, \sum_{j=1}^{n} k_{r,j}= \nu_{r}$, and since there are $r \in K$ number of rows, these can be, in the most simple case imposed using nested loops say $for(\quad;\quad ; )$ loop, with total complexity being the product of complexities of all loops. Each loop runs with a complexity of $O(N)$, where $N$ is number of times of iterations a loop runs. From $\forall r, \sum_{j=1}^{n} k_{r,j}= \nu_{r}$, each $N_{r} = \binom{\nu_{r} + n -1}{n-1}$. This is the formula which counts the solutions of $\sum_{j=1}^{n} k_{r,j}= \nu_{r}$, which can be viewed as placing $n$ identical dots in $\nu_{r}$ number of cells with no restriction to occupancy number of each cell. Therefore, from this perspective,   
$O(\prod_{r \in K} N_{r}) = O( \prod_{r \in K} \binom{\nu_{r} + n -1}{n-1}) $ is the complexity which arises from nested loops. The other factor comes from evaluating the determinant, for which the complexity is in the worst case $O(n^{3})$. For a determinant of Vandermonde the complexity is $O(n^{2})$, $n$ is the determinant order. The equation \ref{eq:19} requires to evaluate mini-vendermonde determents of order $\leq (n-r)$. we shall consider, however the worst case which gives the complexity $ O\left( n^{3}\prod_{r \in K} \binom{\nu_{r} + n -1}{n-1} \right)$. 
\section{Conclusion}
\noindent The Wedge $\wedge$ product is used to find the form of expression for the product of a symmetric $p_{(l)}$ and antisymmetric $V(x_{i})$, which is the basic entity in Frobenius character Formula and is expressed as a sum of scales or determinants which although are alternating polynomials contain or may be made to contain two sub blocks,  one is governed by non-negative integers with sum conditions imposed on each row as $\sum_{\forall r\in K} k_{r,j} = \nu_{r}$, where $k_{r,j} \in \mathbb{N_{0}}$. and the other purely by geometric entries or coordinates $(x_{1}, x_{2},\cdots, x_{n})$. The resulting equation reproduces correctly the character table of $S_{2}$ and $S_{3}$ groups. \\
\newpage
  
\newpage
  \appendix
 \section{First Few Schur Functions and $\chi^{(\lambda)}_{(l)}$ Tables} for $S_{2}$ and $S_{3}$ \label{app:01}
\begin{center}\textbf{First Few Schur Functions}\end{center}
\begin{longtable}{p{0.20\textwidth}|p{0.10\textwidth}|p{0.80\textwidth}}
\caption{$s_{(\lambda)}(x_{1},x_{2},x_{3})$}
\label{t:a1}\\
Young Diagram&$(\lambda)$ & $s_{(\lambda)}$ \\
\hline
$\varPhi$ & (0) & 1\\
\yng(1) &(1)&$x_{1} + x_{2} + x_{3}$			   \\
\yng(2) &(2)&$x_{1}^{2} + x_{1} x_{2} + x_{1} x_{3} + x_{2}^{2} + x_{2} x_{3} + x_{3}^{2}$			   \\
\yng(1,1) &(1,1) &$x_{1} x_{2} + x_{1} x_{3} + x_{2} x_{3}$	  		  \\
\yng(3) &(3) &$x_{1}^{3} + x_{1}^{2} x_{2} + x_{1}^{2} x_{3} + x_{1} x_{2}^{2} + x_{1} x_{2} x_{3} + x_{1} x_{3}^{2} + x_{2}^{3} + x_{2}^{2} x_{3} + x_{2} x_{3}^{2} + x_{3}^{3}$	  		  \\
\yng(2,1) &(2,1) &$x_{1}^{2} x_{2} + x_{1}^{2} x_{3} + x_{1} x_{2}^{2} + 2 x_{1} x_{2} x_{3} + x_{1} x_{3}^{2} + x_{2}^{2} x_{3} + x_{2} x_{3}^{2}$	 		   \\
\yng(1,1,1) &(1,1,1) &$x_{1} x_{2} x_{3}$		    \\
 \hline			
\end{longtable}
\begin{center} \label{c:ps}
\textbf{First Few Power Symmetric Functions in terms of Schur basis}(\ref{eq:15})
\begin{center}
\begin{enumerate}
\item $p_{(1)} = s_{(1)}$
\item $p_{(2)} = s_{(2)} - s_{(1,1)}$
\item $p_{(3)} = s_{(3)} - s_{(2,1)} + s_{(1,1)}$
\end{enumerate}
\end{center}
\end{center}
\begin{longtable}{p{0.10\textwidth}|p{0.10\textwidth}p{0.10\textwidth}}
\caption{$\chi^{(\lambda)}_{(l)}$ for $S_{2}$}
\label{tab:tabs3}\\ \hline
&$(1^{2})$&(2)\\
\hline
\yng(2)&1&1\\
 \yng(1,1)&1&-1\\ 
\hline			
\end{longtable}
\begin{longtable}{p{0.10\textwidth}|p{0.10\textwidth}p{0.10\textwidth}p{0.10\textwidth}}
\caption{$\chi^{(\lambda)}_{(l)}$ for $S_{3}$}
\label{tab:tab}\\ \hline
&$(1^{3})$&(1,2)&(3)\\
\hline
\yng(3)&1&1&1\\
 \yng(2,1)&2&0&-1\\ 
  \yng(1,1,1)&1&-1&1\\ 
\hline			
\end{longtable}
 \section{Calculations for $S_{3}$ Classes} \label{app:02}
\begin{longtable}{p{0.16\textwidth}|p{0.9\textwidth}}
\caption{Table for $(1^{3})$ in $S_{3}$}
\label{tab:t2}\\
\hline
$(k_{1,1}, k_{1,2}, k_{1,3})$& $p_{(1^{3})} V(x_{1}, x_{2},x_{3})$ \\
\hline
(0,0,3)  &  $ p_{(1)}( -x_1^2x_2x_3^2 + x_1x_2^2x_3^2 )$\\
(0,1,2)   &  $ p_{(1)}(  -2x_1^2x_2^2x_3 + x_1^2x_3^3 + 2x_1x_2^3x_3 - x_1x_3^4)$\\
(0,2,1)  &   $ p_{(1)}(  -x_1^2x_2^3 + 2x_1^2x_2x_3^2 + x_1x_2^4 - 2x_1x_2x_3^3)$\\
(0,3,0)  &   $p_{(1)}(   x_1^2x_2^2x_3 - x_1x_2^2x_3^2)$\\
(1,0,2)   &  $p_{(1)}(   -2x_1^3x_2x_3 + 2x_1^2x_2^2x_3 - x_2^2x_3^3 + x_2x_3^4)$\\
(1,1,1)  &   $ p_{(1)}(  -2x_1^3x_2^2 + 2x_1^3x_3^2 + 2x_1^2x_2^3 - 2x_1^2x_3^3 - 2x_2^3x_3^2 + 2x_2^2x _3^3)$\\
(1,2,0)  &    $p_{(1)}(  2x_1^3x_2x_3 - 2x_1^2x_2x_3^2 - x_2^4x_3 + x_2^3x_3^2)$\\
(2,0,1)   &  $ p_{(1)}(  -x_1^4x_2 + x_1^3x_2^2 - 2x_1x_2^2x_3^2 + 2x_1x_2x_3^3)$\\
(2,1,0)   &   $p_{(1)}(  x_1^4x_3 - x_1^3x_3^2 - 2x_1x_2^3x_3 + 2x_1x_2^2x_3^2)$\\
(3,0,0)   &  $ p_{(1)}(  -x_1^2x_2^2x_3 + x_1^2x_2x_3^2)$\\
\hline
&$p_{(1)}(-x_1^4x_2 + x_1^4x_3 - x_1^3x_2^2 + x_1^3x_3^2 + x_1^2x_2^3 - x_1^2x_3^3 + x_1 x_2^4 - x_1x_3^4 - x_2^4x_3 - x_2^3x_3^2 + x_2^2x_3^3 + x_2x_3^4)$\\
& $p_{(1)}(s_{(2)} +s_{(1,1)}) D(x_{1},x_{2},x_{3})$\\
&$(s_{(3)} +2s_{(2,1)} + s_{(1,1,1)}) V(x_{1},x_{2},x_{3})$\\
\end{longtable}
\begin{longtable}{p{0.20\textwidth}|p{0.20\textwidth}|p{0.50\textwidth}}
\caption{Table for $(1,2)$ in $S_{3}$}
\label{tab:t3}\\
\hline
$(k_{1,1}, k_{1,2}, k_{1,3})$&$(k_{2,1}, k_{2,2}, k_{2,3})$& $p_{(1)}p_{(2)} V(x_{1}, x_{2},x_{3})$ \\	
\hline
(0,0,1)  &    (0,0,1)   &    $0$\\
(0,0,1)  &     (0,1,0)   &    $p_{(1)}p_{(2)} (-x_1^2x_2)$\\
(0,0,1)  &     (1,0,0)  &     $p_{(1)}p_{(2)} (x_1x_2^2)$\\
(0,1,0)   &    (0,0,1)  &     $p_{(1)}p_{(2)} (x_1^2x_3$)\\
(0,1,0)   &    (0,1,0)  &     $0$\\
(0,1,0)  &     (1,0,0)  &     $p_{(1)}p_{(2)} (-x_1x_3^2)$\\
(1,0,0)  &     (0,0,1)  &     $p_{(1)}p_{(2)} (-x_2^2x_3)$\\
(1,0,0)   &    (0,1,0)  &     $p_{(1)}p_{(2)} (x_2x_3^2)$\\
(1,0,0)  &     (1,0,0)  &     $0$		\\
\hline
&&$p_{(1)}p_{(2)} (-x_1^2x_2 + x_1^2x_3 + x_1x_2^2 - x_1x_3^2 - x_2^2x_3 + x_2x_3^2)$\\
&&$p_{(1)}p_{(2)} s_{\phi} D(x_{1},x_{2},x_{3})$\\
&& $p_{(1)} (s_{(2)} - s_{(1,1)}) D(x_{1},x_{2},x_{3})$\\
&& $ (s_{(3)} - s_{(1,1,1)}) V(x_{1},x_{2},x_{3})$\\
\end{longtable}
\begin{longtable}{p{0.16\textwidth}|p{0.90\textwidth}}
\caption{Table for $(3)$ in $S_{3}$}
\label{tab:t4}\\
\hline
$(k_{3,1}, k_{3,2}, k_{3,3})$& $p_{(1^{3})} V(x_{1}, x_{2},x_{3})$ \\
\hline
(0,0,1) &  $p_{(3)}( -x_1x_3^2 + x_2x_3^2)$\\
(0,1,0)  &  $p_{(3)}( x_1x_2^2 - x_2^2x_3)$\\
(1,0,0)   & $p_{(3)}( -x_1^2x_2 + x_1^2x_3)$\\
\hline
&$p_{(3)}( -x_1^2x_2 + x_1^2x_3 + x_1x_2^2 - x_1x_3^2 - x_2^2x_3 + x_2x_3^2)$\\
&$p_{(3)}s_{\phi} D(x_{1},x_{2},x_{3})$\\
&$(s_{(3)} -s_{(2,1)} + s_{(1,1,1)}) V(x_{1},x_{2},x_{3})$\\
\end{longtable}


\begin{thebibliography}{}
\bibitem{ref:p1} Amritanshu Prashad, An Introduction to Schur Functions, arXiv:1802.06073(math), 2018.
\bibitem{ref:B7}  Bruce E. Sagan, The Symmetric Group \textit{Representations, Combinatorial Algorithms and Symmetric functions}, $2^{nd} Edition$, Spinger, 2000.
\bibitem{ref:B12} Dudley E. Littlewood, The Theory of Group Characters \textit{And Matrix Representation of Groups} $2^{nd}$ Edition, OXFORD, 1950.
\bibitem{ref:B4} Morton Hammermesh, Group Theory and its Applications to Physical Problems, Dover, 1962.
\bibitem{ref:B2} Palash B. Pal, A Physicist's Introduction to Algebraic Structures $\textit{Vector spaces, Groups, Topological spaces and more}$, Cambridge University Press, 2019.
\bibitem{ref:B11} Peter Szekeres, A Course in Modern Mathematical Physics, CAMBRIDGE (2004).
\bibitem{ref:B3}  Sadri Hassani, Mathematical Physics, $\textit{A Modern Introduction to Its Foundations}$, Springer, 2001.
\bibitem{ref:B5} William Fulton, Young Tableaux, \textit{With Applications to Repsentation theory and Geometry}, CAMBRIDGE, 1997.
\bibitem{ref:B6}  William Fulton, Joe Harris, Repsentation Theory \textit{A First Course}, Springer, 2004.
\end{thebibliography}
\end{document}